\documentclass[journal]{IEEEtran}

\usepackage[utf8]{inputenc}

\usepackage{graphicx}
\usepackage{subfigure}
\usepackage{pgfplots}
\pgfplotsset{compat=1.8}
\usetikzlibrary{shapes,arrows}
\usepackage[margin=0.5in]{geometry}

\usepackage{amsmath}
\usepackage{amsfonts}
\usepackage{amssymb}
\usepackage{array}
\usepackage{mathabx} % Big times
\usepackage{algorithm}
\usepackage{algorithmic}
\usepackage{mathrsfs}

\usepackage{cite}
\usepackage{url} 
\usepackage{enumerate}
\usepackage{url}
\usepackage{lipsum}
\usepackage{soul}

\usepackage{dsfont}
\usepackage{latexsym}
\usepackage{hhline}
\usepackage{color}
\usepackage{textcomp}
\usepackage{hyperref}
\usepackage{booktabs}

\DeclareMathOperator{\argmin}{argmin}

\begin{document}

\title{ Spectral Unmixing with Multiple Dictionaries }

\author{J\'er\'emy E. Cohen and Nicolas~Gillis % <-this % stops a space
\thanks{J.E. Cohen and N.~Gillis (corresponding author) are with the Department of Mathematics and Operational Research, Facult\'e Polytechnique, Universit\'e de Mons, Rue de Houdain 9, 7000 Mons, Belgium. 
E-mails: $\{$jeremy.cohen,nicolas.gillis$\}$@umons.ac.be. The authors acknowledge the support by the F.R.S.-FNRS (incentive grant for scientific research no F.4501.16). NG also acknowledges the support by the ERC (starting grant no 679515). }    
 }

\markboth{}
{}

\maketitle

\begin{abstract} Spectral unmixing aims at recovering the spectral signatures
    of materials, called endmembers, mixed in a hyperspectral or
    multispectral image, along with their abundances. A typical assumption is
    that the image contains one pure pixel per endmember, in which case
    spectral unmixing reduces to identifying these pixels. Many fully automated
    methods have been proposed in recent years, but little work has been done
    to allow users to select areas where pure pixels are present manually
    or using a segmentation algorithm. 
    Additionally, in a non-blind approach, several spectral libraries may be
    available rather than a single one, with a fixed number (or an upper
    or lower bound) of endmembers to
    chose from each. In this paper, we propose a multiple-dictionary
    constrained low-rank matrix approximation model that address these two
    problems. We propose an algorithm to compute this model, dubbed M2PALS, and
    its performance is discussed on both synthetic and real hyperspectral
    images.  \end{abstract}

\begin{IEEEkeywords}
spectral unmixing, 
hyperspectral imaging, 
nonnegative matrix factorization, 
endmember extraction algorithms. 
\end{IEEEkeywords}

%***********************************************************************************************************************************************

\section{Introduction} 
Developments of the remote sensing technology have given
birth to efficient hyperspectral sensors producing high spectral resolution
images, called hyperspectral images (HSIs). Researchers now have at their
disposal a large quantity of HSIs of very different landscapes on the surface of
the Earth, and even of the surface of Mars through the CRISM mission. One of
the several uses for these HSIs is to determine and monitor the chemical
composition of the scene being studied. Indeed, each material has a
recognizable spectral signature that is captured by high resolution spectral
sensors. Various types of similar materials can be identified, which allows,
for instance, the monitoring of the evolution of the deforestation in the
Amazonian forest~\cite{asner2015quantifying}, or of the melting of Alpine snow~\cite{dozier2004multispectral}.

However, because of the high spectral resolution requirements, HSIs lack in
spatial resolution. Therefore, there may be several materials in a single
pixel, resulting in a mixture of their spectral signatures. Spectral unmixing
is interested in finding these spectral signatures given the HSI\@. 
Spectral unmixing has been a central topic for signal
processing research over the last ten to twenty years; 
see, e.g., the surveys~\cite{Jose12,zare2014endmember,Ma14}. 
In
particular, two approaches to spectral unmixing can be distinguished: the blind
approach where no \textit{a priori} information is available other than known
properties of spectra and abundances (such as non-negativity), and the non-blind
approach where spectral libraries containing reference spectra are available. A
list of well-known algorithms for both approaches is given in
Section~\ref{sec:1}.

\paragraph{Contributions} In this paper, we introduce a new model based on
sparse coding that lies between blind and non-blind
approaches,  
designed to offer new possibilities
for dealing with several spectral libraries for a single HSI\@. For instance,
it enables a simple user-control interface for choosing pure pixels areas in an
HSI, based either on manual selection, or semi-supervised or unsupervised segmentation. The proposed
model also allows for using multiple spectral libraries in the non-blind approach. Spectral
signatures are then chosen from a set of known libraries, with the number of atoms to be picked in each dictionary (or an upper or lower bound). 

\paragraph{Outline}
Section~\ref{sec:1} introduces the basic spectral unmixing problem cast in terms of
sparse coding, along with some state-of-the-art methods to tackle it. 
Section~\ref{sec:2} introduces the proposed multiple dictionary matrix
factorization model and a fast algorithm dubbed M2PALS to tackle it featuring an assignment
problem.
Section~\ref{sec:3} shows the efficiency of the proposed model and algorithm on
synthetic data built from the Urban HSI compared to state-of-the-art
methods, and an illustration on the Urban HSI where pure-pixel areas
are either manually selected or automatically 
selected using off-the-shelf segmentation methods.  

\paragraph{Notations}
Matrices are represented by case letters $M$, index sets by calligraphic
letters $\mathcal{K}$, the $i$th column of matrix $M$ is denoted $M(:,i)$
and $M(:,\mathcal{K})$ denotes the subset of columns indexed by $\mathcal{K}$. 
We also denote $\#\mathcal{K}$ the cardinality of $\mathcal{K}$, and $M^T$ the
transpose of $M$. 

\section{A combinatorial formulation of sparse coding}\label{sec:1}
First, before introducing the new multiple dictionary matrix factorization
model, let us recall the single dictionary
formulation.  Let $M$ be a $m$-by-$n$ data matrix. 
In this work, %the following assumptions are imposed on $M$: 
we assume the following: 
\begin{enumerate}
    \item The matrix $M$ admits an approximate low-rank factorization $M \approx AB^T$ of size $r$, that is, $A$ and $B$ have $r$ columns. 
    \item The noise is Gaussian. Missing data, stripes or impulsive noise are
        ignored although common in hyperspectral imaging. The factorization
        therefore can be written as $X = AB^T + N$ where $N$ is the realization
        of a random variable following a white Gaussian distribution. 
    \item Columns of factor $A$ are a subset of columns of the known $m$-by-$d$ dictionary matrix $D$.
\end{enumerate}

These assumptions leads to the following low-rank factorization model for $M$:
%\begin{equation}
%    \begin{array}{cl}
%        M & = AB^T + N \\
%        A & = D\left(:,\mathcal{K}\right) = DS \text{ where } S \in  {\{0,1\}}^{d\times r} \\
%        \text{vec}(N) & \sim \mathcal{N}\left(0, \sigma^2 I_{mn} \right)
%    \end{array}\label{eq:model}
%\end{equation}
\begin{align}
     M & = AB^T + N, \text{vec}(N) \sim \mathcal{N}\left(0, \sigma^2 I_{nm}  \right), \nonumber \\ 
     A & = D\left(:,\mathcal{K}\right) = DS \; \text{ where } S \in  {\{0,1\}}^{d\times r},   \label{eq:model} %\\
     %, \nonumber  
\end{align}
for a given noise variance $\sigma^2$. In this model, $\mathcal{K}$ is a set of
indices of atoms in $D$ corresponding to the columns of the factor matrix $A$. 
The matrix $S$ is a sparse selection matrix which has only 1 non-zero entry per
column. 

The literature on computing $S$ and $B$ is extensive. A first group of methods
are the continuous methods, based on applying iterative descent algorithms like
ADMM~\cite{ammanouil2014glup} or fast gradient~\cite{GL16} to a relaxed version
of~\eqref{eq:model}. In most works, the product of $S$ and $B$ is estimated directly, 
using another variable $X=SB^T$ that has to be row sparse, with $M \approx DX$.  
This constraint on $X$ can be relaxed
using the $\ell_1$ norm, the mixed norm $\ell_{2}/\ell_{1}$ or other convex
norms, making the problem convex~\cite{EMO12, elhamifar2012see, BRRT12, GL16}.  
However, it has a very large number of parameters to be
optimized ($X$ is a $d$-by-$n$ matrix) so that continuous methods may be slow and memory consuming. 
Note that continuous methods allow for adding various constraints like non-negativity.

A second category of sparse coding methods are greedy, or not based on a
continuous optimization algorithm:  
\begin{itemize}

    \item A large family of %non-iterative 
    geometric algorithms have been
developed in the case where $D$ is the data $M$ itself, which is a common
assumption in spectral unmixing referred to as ``the pure-pixel assumption''.
This assumptions implies that the spectra in $M$ form a simplex with $r$ vertices. 
%, so that spectral unmixing reduces to identifying these vetices. 
These geometric algorithms include vertex component analysis (VCA), successive projection algorithm (SPA), 
successive nonnegative projection algorithm (SNPA) and NFINDR, to cite a few, which aim at finding these
vertices by resorting to geometric tools, such as
projections~\cite{Win99,nascimento2005vertex,G14b}. 
%Non-negativity is mandatory in SNPA and \ldots.
These methods are usually very fast even for very large data set, but are
restricted to the pure-pixel case and may not behave well when the data is grossly corrupted. 
Also the $\ell_2$ reconstruction error may be high since no least
square criterion is minimized directly. 

    \item In the general case where the dictionary is not necessarily
        the data itself, one of the many efficient methods for solving
        the general sparse coding problem is the Simultaneous Orthogonal
        Matching Pursuit
        algorithm~\cite{tropp2006algorithms}. It has been referred to in the
        spectral unmixing community as Self Dictionary SOMP~\cite{Fu15}.
        Recently, Matching Pursuit Alternating Least Squares has
        been proposed that aims at minimizing a $\ell_2$ cost function while
        using a fast greedy approach to find the right atoms in the
        dictionary~\cite{cohen2017new}. 

\end{itemize}

Finally, a last family of methods specific to the spectral unmixing problem
uses smart variations on the naive brute force algorithm that computes all
possible combinations of atoms to be selected for each pixel. These methods
include MESMA, MESLUM, AUTOMCU, or more recently,
AMUSES~\cite{roberts1998mapping,combe2008analysis,asner2003scale,degerickx2017novel}. A direct consequence of the pixel-wise
approach is that the low-rank hypothesis may be
violated. % *** je ne comprends pas trop cette phrase? 
Also when the spectral library is large, such methods may be
time-consuming. On the other hand, algorithms like MESMA may handle large
scenes featuring more than a hundred endmembers. Because these methods may not
produce results that satisfy our required assumptions, they will not be used
for comparison in this short paper.

In the next section, we introduce another model closely related
to~\eqref{eq:model} which allows multiple dictionaries to be used at once as a set of
admissible spectral signatures.

\section{Multiple-Dictionary Matrix Factorization}\label{sec:2}

In the self-dictionary setting, a user may want some control over the regions
of the HSI where the pure pixels are selected from. On the other hand,
hand-picking pixels that seem pure may lead to poor results because the
pure-pixel property is difficult to assess visually. Moreover, using a
segmentation algorithm to compute homogeneous areas does not provide a good
input for usual sparse coding method since the segmented regions would be bundled,
and coherence lost. Similarly, when working
with external libraries, it is reasonable to assume that only a few materials
from a specific library should be selected. For example, exactly, at
least or at most two spectra
related to vegetation among some 20 available vegetation spectra. Lumping the
libraries together and running any methods described above %(except the brute-force family) 
would not guarantee such an assignment.

To tackle these issues, we suggest to drop the single dictionary constraint.
Using notations from~\eqref{eq:model}, the third working hypothesis is modified
to:
\begin{itemize}
    \item[3.] Columns of matrix $A$ belong to matrices $D_i$, with
        ${d_i}$ the exact number of atoms to be picked in each $D_i$, $1 \leq i \leq p$.
\end{itemize} 
This new multi-dictionary model becomes: 
%\begin{equation}
 %   \begin{array}{cl}
 %       M & = AB^T + N \\
 %       A & =
 %       \left[D_{u_1}\left(:,\mathcal{K}_{u_1}\right),\dots,D_{u_r}\left(:,\mathcal{K}_{u_r}\right)\right] \\
 %       \# \mathcal{K}_{u_i} & = d_{u_i} \text{ and } \sum\limits_{i=1}^{r}{d_{u_i}} = r \\
 %       \text{vec}(N) & \sim \mathcal{N}\left(0, \sigma I_{nm} \right)
 %   \end{array}\label{eq:model2}
%\end{equation}
\begin{align}
      M & = AB^T + N, \, \text{vec}(N) \sim \mathcal{N}\left(0, \sigma I_{nm} \right), \nonumber \\
     A & =
     \left[D_{1}\left(:,\mathcal{K}_{1}\right),\dots,D_{p}\left(:,\mathcal{K}_{p}\right)\right]\Pi, \label{eq:model2} \\
     \# \mathcal{K}_{i} & = d_i \text{ (or } \leq d_i \text{) and } \sum\limits_{i=1}^{p}{\#
      \mathcal{K}i}
     = r,  \nonumber %\\
       %\text{vec}(N) & \sim \mathcal{N}\left(0, \sigma I_{nm} \right),  \nonumber 
\end{align}
%where $u_i$ is the index in $[1,p]$ of the dictionary from which $A(:,i)$ ispicked, and 
where $\mathcal{K}_i$ contains the indices of the $d_i$ selected atoms in
dictionary $D_i$ for $1 \leq i \leq p$, and $\Pi$ is the permutation matrix
that matches factors to their corresponding atoms. 
Note that the model~\eqref{eq:model2} can be adapted to handle different
situations; in particular if we have a lower bound  on the number of atoms to
be selected from each dictionary (instead of the exact value or an upper
bound).  
In that case, we introduce a new dictionary $D_{p+1}$ containing all the other ones and choose $d_{p+1} = r - \sum_{i=1}^p d_i$ 
where $d_i$ is the lower bound for the $i$th dictionary.  
Let us illustrate this on a simple example: assume we are given two dictionaries $D_1$ and $D_2$ from which we have to select three endmembers. 
We only know that at least one endmember has to be picked from each dictionary (that is, $d_1 \geq 1$ and $d_2 \geq 1$). 
In that case, we introduce a third dictionary $D_3 = [D_1, D_2]$ and impose $d_1=d_2=d_3=1$ in model~\eqref{eq:model2}. 

% *** je n'ai pas compris pq tu as introduit le u_i, au lieu de prendre juste _i comme dans la phrasejuste au dessus
% De plus il y avait une erreur: r au lieu de p 

To compute the sets $\mathcal{K}_{i}$'s and the abundance matrix $B$, we choose to adapt
the MPALS algorithm introduced in~\cite{cohen2017new}. The reason is two-fold: 
(1)~MPALS was shown to perform extremely well, outperforming geometric and continuous approaches in many cases, and 
(2)~it can be easily adapted to the multiple dictionary scenario by solving an assignment problem (see below).  
The proposed algorithm,
named Multiple Matching Pursuit Alternating Least Squares (M2PALS), is an alternating algorithm that computes the second factor $B$
using a least squares update, while using a least squares estimate of $A$ to
serve as a proxy for computing the sets and indices $\mathcal{K}_{i}$. 
It aims at solving 
\[ 
    \underset{A=[D_{1}(:,\mathcal{K}_{1}),\dots,D_{p}(:,\mathcal{K}_{p})],~B}{\text{argmin}}   
    \quad \| M - AB^T \|_F, 
\] 
by first estimating $A$ and $B$ without the dictionary constraint, then projecting $A$ onto the atoms of the dictionaries. 
As we will see, the projection is a standard assignment problem. 

\begin{algorithm}
    \begin{algorithmic}
        \STATE{\textbf{Input: }Data matrix $M$, dictionaries $D_1,\dots,D_p$,
            exact (or upper bound on the) number of atoms to pick in
        each library $d_1,\dots,d_p$, initials factors $A$ and $B$.}
        \STATE{\textbf{Output: }Factors $A$ and $B$, selected atoms
        $\mathcal{K}_{1},\dots,\mathcal{K}_{p}$. } %($1 \leq i \leq p$).} 
        \WHILE{the convergence criterion is not met}
            \STATE{\ul{Least squares estimate of $A$}: $A = \argmin_{X} ||M - X B^T||_F$}  % $A = MB^\dagger$ **** Je n'aime pas utiliser dagger car (1) ce n'est pas très explicite, (2) en pratique, il ne faut (surtout) pas calculer le pseudoinverse pour résoudre ce problème numériquement
            \STATE{\ul{Solve the assignment problem}: Find the $\mathcal{K}_{i}$'s to match $A$ the best using the available dictionaries, that is, 
                solve~\eqref{eq:assignement}, and update $A
            =[D_{1}(:,\mathcal{K}_{1}),\dots,D_{p}(:,\mathcal{K}_{p})]\Pi$} 
            \STATE{\ul{Least squares estimate of $B$}: $B = \argmin_{Y} ||M - A Y^T||_F$ } % $B = A^\dagger M$ 
        \ENDWHILE\end{algorithmic}
    \caption{M2PALS} 
\end{algorithm}
If non-negativity is imposed on $A$ and/or $B$, M2PALS can be adapted into
M2PNALS by using a non-negative least squares solver to update $A$ and
$B$. Typically, these will not be solved up to global optimality for faster
computation. In this paper, we will use a block coordinate descent method;
see~\cite{gillis2012accelerated}. Other constraints on the abundance maps, for
instance, homogeneity constraints using total variation
regularization~\cite{iordache2012total}, can
be easily added in M2PALS since the constrained estimation of abundances $B$
knowing $A$ is simply a
constrained least squares problem.

\paragraph{The assignment problem} 

In M2PALS, after $A$ has been estimated through a least squares update, one has
to compute the index sets $\mathcal{K}_{i}$'s, that is, for each column of $A$, 
find the nearest column from one
of the dictionaries $D_i$, while satisfying the constraints that exactly
(or at most) 
$d_i$ columns can be selected in the dictionary $D_i$. 
If only one dictionary is
provided, each column of $A$ can be processed independently, by picking the nearest column according to some distance criterion.  
However with multiple dictionaries and the constraint on the fixed number of atoms to be selected in each dictionary,  
the following assignment problem should be solved: 
Given $A \in \mathbb{R}^{m \times n}$, 
find the index sets $\mathcal{K}_i$ ($1 \leq i \leq p$) with $\#\mathcal{K}_i = d_i$ and $\sum_{i=1}^p d_i = r$,  
and find a permutation matrix $\Pi \in {\{0,1\}}^{r \times r}$ such that the matrix $\tilde{A} = [D_{1}(:,\mathcal{K}_{1}),\dots,D_{p}(:,\mathcal{K}_{p})] \Pi$ minimizes 
\begin{equation}
  \sum_{j=1}^n \text{dist} \left( A(:,j),  \tilde{A}(:,j)  \right), 
\label{eq:assignement}
\end{equation}
 where dist$(.,.)$ is some distance criterion. Note that the permutation matrix $\Pi$ is necessary since the columns of the input matrix $A$ are not necessarily ordered in the same way as the dictionaries. 
 The problem can be solved as follows. 
 If we assume we know the dictionary from which we pick the atom to approximate each column of $A$, then the problem is trivial: for each column of $A$, we pick the column of the given dictionary which is the closest according to dist$(.,.)$. 
 It remains to identify which dictionary is used to approximate each column of $A$. 
 Defining the distance between the $j$th column of $A$ and the $i$th dictionary as 
 \[ 
 d(i,j) = \min_{k} \; \text{dist} \left( A(:,j) , D_i(:,k) \right), 
 \] 
 problem~\eqref{eq:assignement} reduces to an assignment problem where 
 on the left hand side we have the $r$ columns of $A$, 
 on the right hand side we have $d_i$ copies of dictionary $D_i$ for a total of $r$ nodes (recall $\sum_{i=1}^p d_i = r$), 
 and the weights are given by $d(i,j)$. If $d_i$ are upper bounds on the
 number of atoms to be selected, then the total number of nodes is larger than
 $r$ but the assignment problem is solved in a similar fashion.
Indeed, recall that an assignment problem is defined as follows: 
 given $r$ persons, $l$ tasks ($l \leq r$) and the cost $c_{i,j}$ to assign the $i$th person to the $j$th task, 
assign a different person to each task in order to minimize the total cost. Mathematically, it can be formulated as follows 
 \begin{equation*}
 \min_{x_{i,j} \in \{0,1\}, 1 \leq i,j \leq r}  \; \sum_{i,j} c_{i,j} x_{i,j} \; 
 \text{ s.t. }  \;     \sum_i x_{ij} \leq 1 
 \text{ and } \sum_j
 x_{ij} = 1, 
 \end{equation*}
 where $x_{ij} = 1$ if and only if person $i$ is assigned to task $j$. 
 An optimal assignment can be obtained using linear programming, or the
 Hungarian algorithm\footnote{We use the code available at\\ \url{http://www.mathworks.com/matlabcentral/fileexchange/20328-munkres-assignment-algorithm}.}~\cite{kuhn1955hungarian}. 

Several reasonable distances can be used in the context of spectral unmixing.
The Euclidean distance is a simple choice but may be sensitive to scaling of spectra (that can occur for example because of inhomogeneous light intensity within the HSI). 
Many authors prefer to resort to an angle measurement
through the SAM distance, 
or to remove the average of the spectra before computing
an angle measurement, which yields the mean-removed spectral angle (MRSA) distance. 
%Given two spectral signatures, $x, y \in \mathbb{R}^m$, the MRSA is defined as 
%\begin{equation} \nonumber 
%\phi(x,y) 
%= \frac{100}{\pi}
%\arccos \left( \frac{ {(x-\bar{x})}^T (y-\bar{y}) }{||x-\bar{x}||_2 ||y-\bar{y}||_2} \right) \quad \in \quad [0,100],  
%\end{equation}
%where, for a vector $z \in \mathbb{R}^m$,  $\bar{z} = \frac{1}{m} \left(\sum_{i=1}^m z(i)\right) e$ with $z(i)$ being the $i$th entry of $z$ and $e$ the vector of all ones. The MRSA measures how close two endmembers are (neglecting scaling and translation); 0 meaning that they match perfectly, 100 that they do not match at all.  
These metrics are equally fast to compute, thus the choice is entirely application dependent.
Note that various distance measurements may result in various assignment
results; see, e.g., the discussion in~\cite{keshava2004distance}.  
In the simulation section, we will use the normalized inner product to asses the closeness between two vectors, 
that is, $\text{dist}(x,y) = 1-\frac{x^T y}{||x||_2 ||y||_2}$, 
which is equivalent to maximize the sum of the normalized inner products
between the selected atoms of the dictionaries and the columns of $A$. 
If the atoms of the dictionaries are normalized to have unit $\ell_2$ norm, it can
be shown that the normalized inner product, plugged in the assignment problem, provides the
best solution for the index sets $\mathcal{K}_i$ for matching $A$ in the least
square sense. 

%\begin{remark} The multiple dictionary low-rank factorization may be seen as a more
%constrained version of the usual sparse coding problem. Indeed, a more relaxed
%problem is obtained if all dictionaries are stacked into one large library.
%Therefore, results on identifiability and existence of the best approximate
%solution obtained in~\cite{cohen2017new} propagate to the multi-dictionary model.
%\end{remark}

\section{Application for semi-supervised and unsupervised pure-pixel selection}\label{sec:3}

This section is divided in two parts. 
The first one is a comparison of M2PNALS (the
non-negative version of M2PALS) on synthetic data with state-of-the-art spectral unmixing algorithms,
namely 
SPA~\cite{MC01}, 
SNPA~\cite{G14b}, 
NFINDR~\cite{Win99}, 
Fast Gradient for Nonnegative Sparse Regression (FGNSR)~\cite{GL16}, 
Group Lasso Unit sum Positivity constraints~\cite{ammanouil2014glup}, 
SDSOMP~\cite{Fu15} and
MPALS~\cite{cohen2017new}. 
A recent segmentation algorithm, hierarchical rank-one NMF (H2NMF)~\cite{gillis2015hierarchical}, is also used as a preprocessing
tool for M2PNALS to compute a set of dictionaries in an unsupervised manner.

The second one is an example of using M2PNALS to develop a user-friendly pure-pixel selection
interface, based on either manual or automatic selection of image areas
containing pure pixels.

\subsection{Experiments on synthetic HSIs} 

For this experiment, we chose $r=6$
spectra from the Urban HSI\footnote{Available on \url{http://www.erdc.usace.army.mil/}} (162 spectral bands, images 307 by 307 pixels), which were selected using SPA on the
original image. Then the following setting is repeated $N=100$ times:
abundances are drawn from a uniform Dirichlet distribution so that they belong to a
simplex with $r$ vertices, and the $r$ pure pixels are added at random positions in the data set.  
The data is generated by adding a white Gaussian noise with fixed signal to noise
ratio (SNR), and then clipping negative entries to zero\footnote{Clipping to zero implies that the
SNR is slightly over-estimated, as the noise is no longer Gaussian.}. 
The size of the abundance matrix is set to
$n=200$, which is rather small but allows all methods to run in a
reasonable computational time. 
All tests are preformed using Matlab on a laptop Intel CORE i5-3210M CPU
@2.5GHz 6GB RAM\@. 

At various SNR levels, for each setting and for each algorithm, the number of
wrongly-selected pixels in the synthetic image is stored. 
Average results are presented in Table~\ref{table:sim} in percentage of miss-selected atoms, 
along with the average run times.  
For M2PNALS, additionally, the size of the dictionaries is used as a grid parameter. For
instance, M2PNALS-10 means that each dictionary contains 1 pure-pixel plus 9
other randomly chosen spectra from the synthetic HSI (which may be one of the
other pure pixels). The M2PNALS-01 is provided for sanity check and is
trivially perfect for selecting pure pixels since each dictionary contains
exactly one pure pixel. 

Not too surprisingly, M2PNALS performs better than state-of-the-art methods in
the noisy cases. This is natural since M2PNALS makes use of \textit{a priori}
information provided by the user on the location of the pure pixels.  In the
noiseless case, all algorithms perform perfectly except the SDSOMP algorithm,
which is a sparse-coding method not adapted to identifying endmembers and
rather selects the most representative spectra in the data at each iteration.
Also, using H2NMF to segment the synthetic image (here without using any
spatial information) and provide candidates $D_i$ deteriorates the results
at high SNR but improves on identification at medium SNR\@. This is natural since
segmentation here can be understood as a regularization of the pure pixel
selection process. In this simulated experiment, the data lies exactly in a
simplex, thus linear unsupervised methods should perform better at high SNR\@. On the other
hand, the segmentation process may include two pure pixels in the same
segmented cluster, so that we set $d_i=2$ for H2NMF-M2PNALS\@. See the next experiment for a more complete discussion on
using segmentation before M2PNALS\@.  

Moreover, note that MPNALS is an algorithm sensitive to
initialization although it is able to identify good solutions in most cases regardless of the initialization; 
see the discussion and experiments in~\cite{cohen2017new}.  
In this experiment, we chose SNPA for initialization, but
better results may be achieved with a different initialization choice. 
The same comment applies to M2PNALS since it is adapted from MPNALS\@.  Also, it can be
checked that M2PNALS is fast compared to continuous algorithms. Actually, its
complexity per iteration is linear in $m$, $n$ and $r$ (assuming we use a
first-order method to solve the nonnegative least squares subproblems), which
is the same as geometrical algorithms.  Therefore, it can be used on large
images, as illustrated in the next experiment. 

\begin{table}\label{table:sim}
    \centering
    \begin{tabular}{c|ccccc|c}
    SNR & 0 & 10 & 20 & 30 & 50 & Time (s.)   \\
    \toprule
    SPA        & 93.66 & 79.16 & 51.83 & 13.33 & 0 & 0.0009 \\ 
    MPANLS     & 93.50 & 80.16 & 47.83 & 1.66 & 0 & 0.0955 \\ 
    SNPA       & 93.66 & 79.16 & 51.33 & 2.00 & 0 & 0.0420 \\
    FGNSR      & 97.00 & 86.16 & 58.16 & 17.66 & 0&  3.7960 \\ 
    GLUP       & 96.16 & 82.00 & 31.83 & 6.66 & 0 & 1.1057 \\     
    NFINDR     & 96.00 & 85.50 & 43.50 & 1.66 & 0 & 0.0192 \\
    SDSOMP     & 95.83 & 92.66 & 94.50 & 61.50 & 16.83 & 0.0154   \\
    M2PNALS-01 & 0 & 0 & 0 & 0 & 0 & 0.1068 \\  
    M2PNALS-10 & 67.16 & 59.33 & 42.16 & 0.50 & 0 & 0.1161 \\
    M2PNALS-25 & 83.33 & 74.16 & 51.16 & 0.50 & 0 & 0.1240 \\
    M2PNALS-50 & 88.66 & 79.50 & 54.66 & 0.50 & 0 & 0.1332 \\
    H2NMF-M2PNALS & 93.66 & 80.33 & 48.00 &  3.00    & 1.50  & 0.1057 \\
\end{tabular} 
\caption{Percentage of miss-selected atoms in synthetic HSIs depending on the SNR (first five columns), 
and average run times in seconds (the stopping criterion for MPNALS and M2PNALS is that the reconstruction
error varies less than $10^{-3}\%$ in two successive iterations 
which always happened within 50 iterations in this simulation).}  
\end{table}

\subsection{Using M2PALS for a user selection of pure-pixel areas} 

Our motivation for developing M2PNALS is not to outperform state-of-the-art
spectral unmixing methods in all scenarios. Rather, M2PNALS is meant to open new
possibilities for dealing with spectral unmixing of HSI where a user inputs
several spectral libraries.

The flexibility of M2PNALS lies in the variety of ways that the dictionaries can be
defined. Of course they can be given by a user as external spectral libraries,
but at least three other ways come to mind that only involve the available HSI.
Although well-known dictionary learning methods such as~\cite{jiang2011learning} 
do not apply here since only one dictionary would be learned this way, any
segmentation or classification algorithm, which provides a labelled partition of
the pixels, can be used to cut the original image into several regions that can
be used directly as spectral libraries. Below, we introduce three different
strategies for segmentation: a fully supervised, hand selection method, where
the user chooses areas that she or he believes contain pure pixels; a
semi-supervised state-of-the-art segmentation algorithm~\cite{li2012spectral} based on
logistic regression, that makes use
of spatial correlations, and for which the
training data has been labelled by hand; an unsupervised recent algorithm
H2NMF~\cite{gillis2015hierarchical} where no spatial information is used to perform segmentation. 

\begin{figure*} \centering  
    \includegraphics[width=0.32\textwidth]{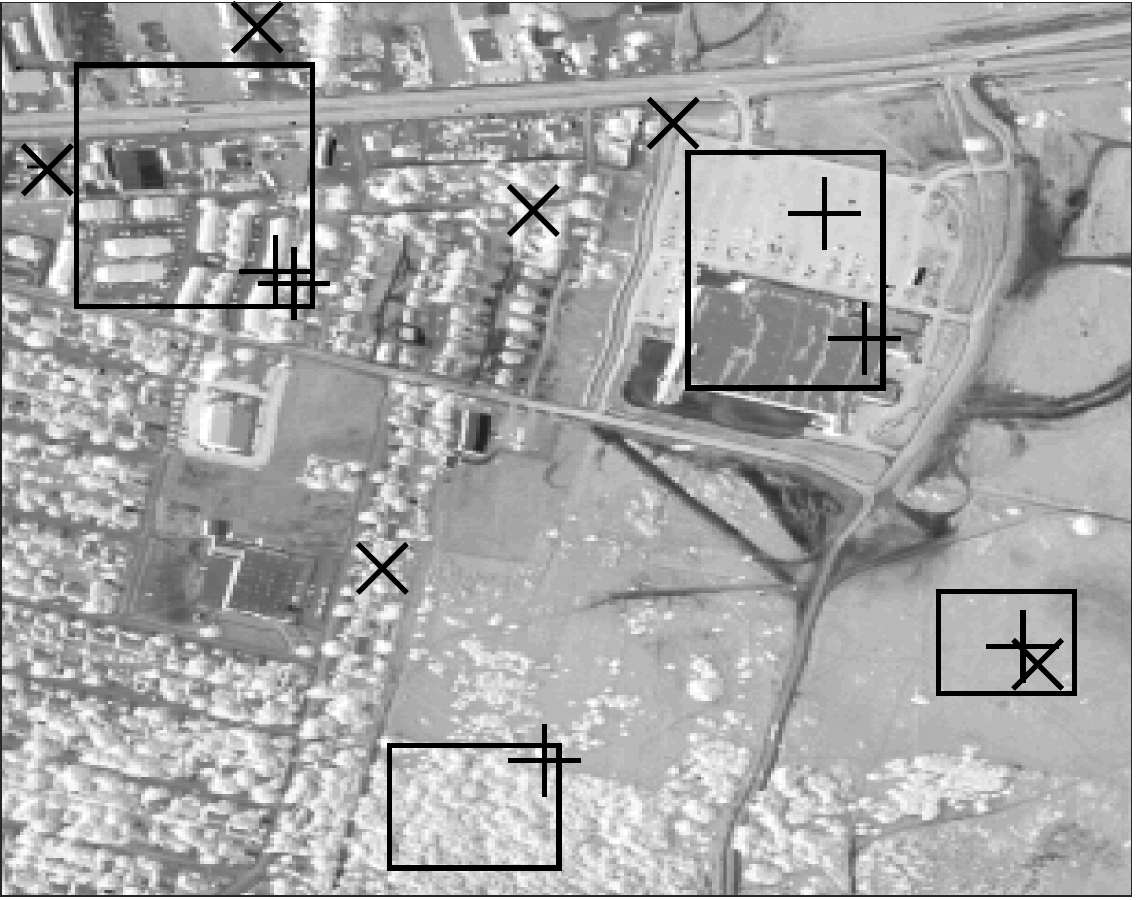}
    \includegraphics[width=0.32\textwidth]{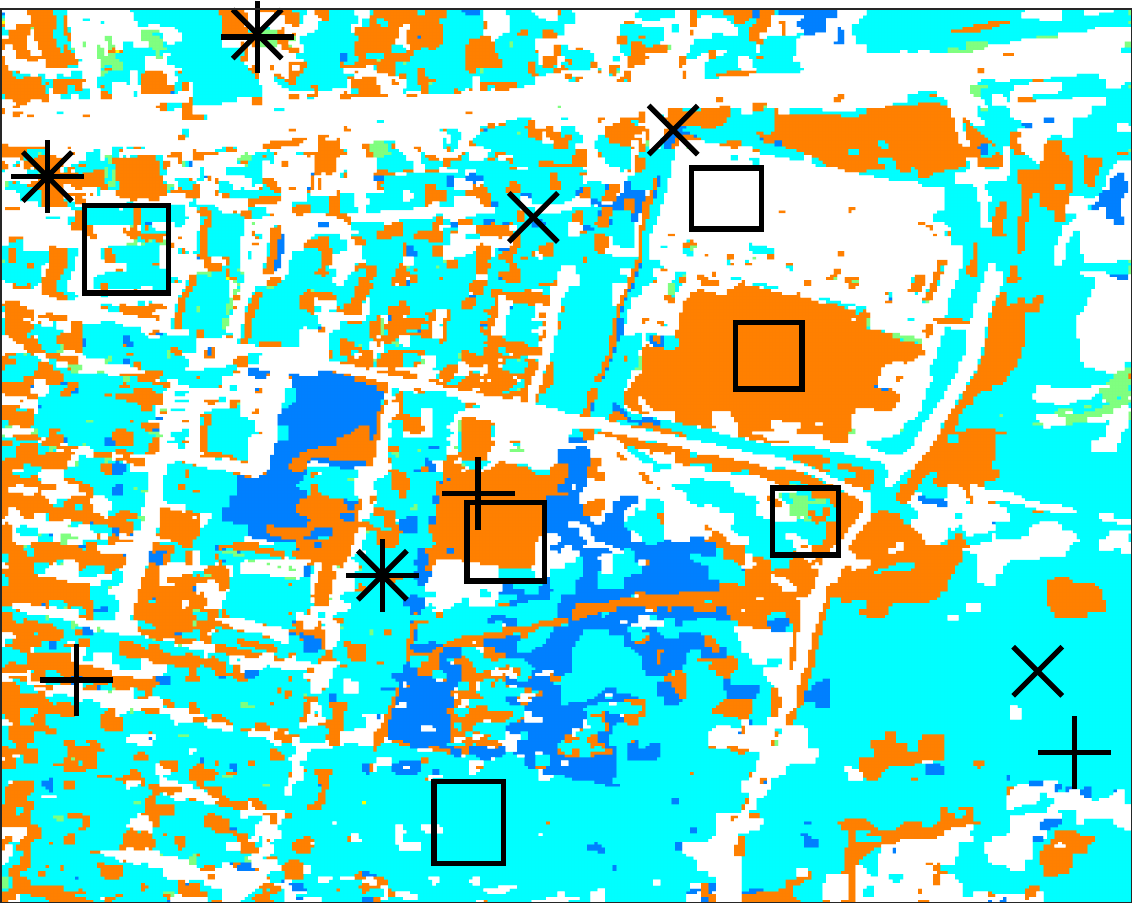}
    \includegraphics[width=0.32\textwidth]{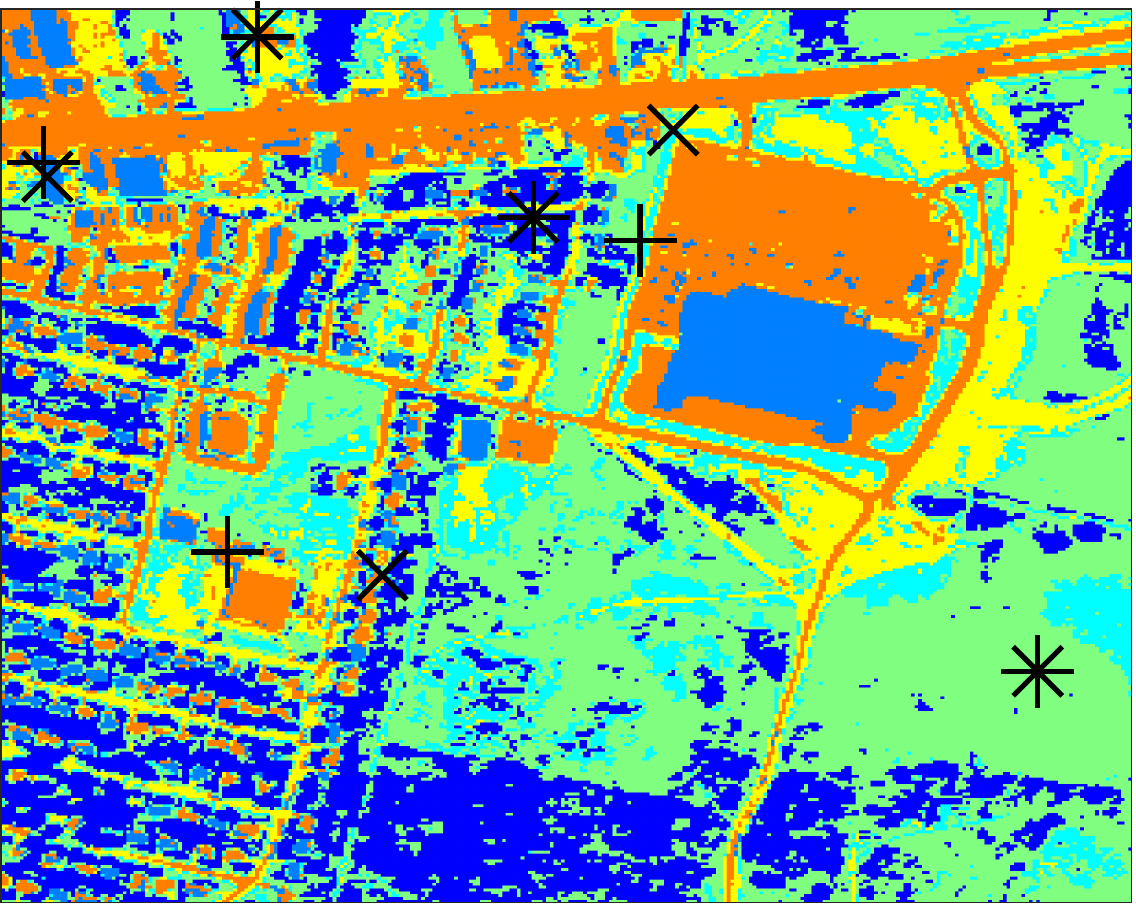}
    \caption{ An example of a manual (left), semi-supervised (middle) and
        unsupervised selection (right) of pure pixel areas in the Urban
    HSI\@. In the manual selection scenario, squares represent
    two areas containing two pure pixels
    and two areas containing one. The squares in the semi-supervised
    segmentation case represent the training data.  
		The pixels identified by M2PNALS are represented with plusses `{\large{+}}' and lie
    within the selected areas, 
		while pixels selected by fully-blind MPALS
    are represented with crosses `{\large{$\times$}}'.  
		}\label{fig:select}
\end{figure*}

 These methods are illustrated  in this experiment on the
Urban HSI using $r=6$ endmembers, along with the MPALS algorithm, see Figure~\ref{fig:select}. For the manual selection method,
 the user is asked to determine the number $d_i$
of pure pixels in each area, whereas $d_i$ are fixed to $1$ for the
unsupervised segmentations, and to $2$ for the semi-supervised one since one of
the cluster contained very few atoms. In this example, the centroids
of the H2NMF clusters were used  to initialize both
MPNALS and M2PNALS\@. The relative reconstruction error of manual,
semi-supervised and unsupervised M2PNALS here is respectively
$4.28\%$, $4.12\%$ and $4.20\%$ , while that of MPNALS is $4.05\%$ for a running time
of respectively $50$, $47s$, $51s$ and $20s$.  
In all cases, M2PNALS provided a sightly worse reconstruction error, because it is more constrained than MPNALS 
(which can choose any pixels in the HSI) while the selected areas do not
contain the best set of pure pixels (as in the synthetic data experiment).
However, segmentation used in conjunction with M2PNALS allowed the user to control in which areas
pure pixels were selected, which is not possible with MPNALS nor any geometric pure pixel selection method. 

% Dm has been stored in Dm_Urban.mat
% ** Note: add MPNALS with $D=[D_1,\dots,D_p]$? **

\section*{Conclusion}

In this paper, we introduced a new sparse coding  model that allows
for  computing spectral unmixing using unknown spectra that belong to
several classes of known spectra. This model may be used in post-processing,
for instance, after selecting areas containing pure pixels, or after segmentation algorithms that provide
homogeneous areas.  
It may also be used for estimating abundances given an HSI
and a collection of libraries, as an alternative to methods such as MESMA\@. 
The M2PALS algorithm has been developed to identify the parameters of this
multiple dictionary matrix factorization model, and its performance have been demonstrated on a
simulated HSI, and on the Urban HSI\@.  

Further research includes dealing with other types of noise, 
and the study of the identifiability 
of the multiple dictionary model~\eqref{eq:model}, which should be stronger than that of the
single dictionary model~\eqref{eq:model2}; see~\cite{cohen2017new}. 
Also, the authors are working on adapting M2PALS for
the fusion of multispectral and hyperspectral images under the pure-pixel
assumption, since the mapping from an HSI to a MSI induces a pixel to area
mapping.

\section*{Acknowledgments} 

The authors wish to thank David Kun for sharing the HyperSpectralToolbox,
Antonio Plaza for sharing the MLR segmentation code  and Rita
Ammanouil for sharing the GLUP algorithm.

\bibliographystyle{IEEEtran}

\end{document}